\documentclass{amsart}

\usepackage{amsmath,amsthm,amsfonts,amscd,amssymb}
\usepackage[pdfborder={0 0 0}]{hyperref}
\usepackage{enumerate,verbatim}
\usepackage{enumitem}
\usepackage{color}

\newcommand{\ZZ}{\mathbb{Z}}

\newtheorem{thm}{Theorem}[section]
\newtheorem{cor}[thm]{Corollary}
\newtheorem{lem}[thm]{Lemma}
\newtheorem{prop}[thm]{Proposition}

\theoremstyle{definition}

\theoremstyle{remark}
\newtheorem{rem}{Remark}[section]

\begin{document}

\title{Counting ideals in polynomial rings}
\author{Lenny Fukshansky}
\author{Stefan K\"uhnlein}
\author{Rebecca Schwerdt}
\thanks{The first author was partially supported by the NSA grant H98230-1510051.}
\address{Department of Mathematics, 850 Columbia Avenue, Claremont McKenna College, Claremont, CA 91711}
\email{lenny@cmc.edu}
\address{Institut f\"ur Algebra und Geometrie, KIT, FRG-76128 Karlsruhe, Germany}
\email{stefan.kuehnlein@kit.edu}
\address{Arbeitsgruppe f\"ur Kryptographie und Sicherheit, KIT-Fakult\"at f\"ur Informatik, KIT, FRG-76128 Karlsruhe, Germany}
\email{rebecca.schwerdt@kit.edu}


\begin{abstract} 
We investigate properties of zeta functions of polynomial rings and their quotients, generalizing and extending some classical results about Dedekind zeta functions of number fields. By an application of Delange's version of the Ikehara Tauberian Theorem, we are then able to determine the asymptotic order of the ideal counting function in such rings. As a result, we produce counting estimates on ideal lattices of bounded determinant coming from fixed number fields, as well as density estimates for any ideal lattices among all sublattices of $\ZZ^d$. We conclude with some more general speculations and open questions.
\end{abstract}

\maketitle

\def\A{{\mathcal A}}
\def\AA{{\mathfrak A}}
\def\B{{\mathcal B}}
\def\C{{\mathcal C}}
\def\D{{\mathcal D}}
\def\EE{{\mathfrak E}}
\def\F{{\mathcal F}}
\def\G{{\mathcal G}}
\def\x{{\mathcal H}}
\def\I{{\mathcal I}}
\def\II{{\mathfrak I}}
\def\J{{\mathcal J}}
\def\K{{\mathcal K}}
\def\kk{{\mathfrak K}}
\def\L{{\mathcal L}}
\def\LL{{\mathfrak L}}
\def\M{{\mathcal M}}
\def\mm{{\mathfrak m}}
\def\MM{{\mathfrak M}}
\def\N{{\mathcal N}}
\def\O{{\mathcal O}}
\def\OO{{\mathfrak O}}
\def\PP{{\mathfrak P}}
\def\R{{\mathcal R}}
\def\PNR{{\mathcal P_N(\real)}}
\def\PMNR{{\mathcal P^M_N(\real)}}
\def\PdNR{{\mathcal P^d_N(\real)}}
\def\s{{\mathcal S}}
\def\V{{\mathcal V}}
\def\X{{\mathcal X}}
\def\Y{{\mathcal Y}}
\def\Z{{\mathcal Z}}
\def\H{{\mathcal H}}
\def\cee{{\mathbb C}}
\def\Nn{{\mathbb N}}
\def\pee{{\mathbb P}}
\def\que{{\mathbb Q}}
\def\QQ{{\mathbb Q}}
\def\real{{\mathbb R}}
\def\RR{{\mathbb R}}
\def\zed{{\mathbb Z}}
\def\ZZ{{\mathbb Z}}
\def\aaa{{\mathbb A}}
\def\ff{{\mathbb F}}
\def\HDelta{{\it \Delta}}
\def\kk{{\mathfrak K}}
\def\qbar{{\overline{\mathbb Q}}}
\def\kbar{{\overline{K}}}
\def\ybar{{\overline{Y}}}
\def\kkbar{{\overline{\mathfrak K}}}
\def\ubar{{\overline{U}}}
\def\eps{{\varepsilon}}
\def\beps{{\boldsymbol \varepsilon}}
\def\ahat{{\hat \alpha}}
\def\bhat{{\hat \beta}}
\def\gt{{\tilde \gamma}}
\def\h{{\tfrac12}}
\def\be{{\boldsymbol e}}
\def\bei{{\boldsymbol e_i}}
\def\bc{{\boldsymbol c}}
\def\bm{{\boldsymbol m}}
\def\bk{{\boldsymbol k}}
\def\bi{{\boldsymbol i}}
\def\bl{{\boldsymbol l}}
\def\bq{{\boldsymbol q}}
\def\bu{{\boldsymbol u}}
\def\bt{{\boldsymbol t}}
\def\bs{{\boldsymbol s}}
\def\bv{{\boldsymbol v}}
\def\bw{{\boldsymbol w}}
\def\bx{{\boldsymbol x}}
\def\bX{{\boldsymbol X}}
\def\bz{{\boldsymbol z}}
\def\bwy{{\boldsymbol y}}
\def\bY{{\boldsymbol Y}}
\def\bL{{\boldsymbol L}}
\def\ba{{\boldsymbol a}}
\def\bb{{\boldsymbol b}}
\def\bet{{\boldsymbol\eta}}
\def\bxi{{\boldsymbol\xi}}
\def\bo{{\boldsymbol 0}}
\def\bone{{\boldsymbol 1}}
\def\bol{{\boldsymbol 1}_L}
\def\ep{\varepsilon}
\def\p{\boldsymbol\varphi}
\def\q{\boldsymbol\psi}
\def\rank{\operatorname{rank}}
\def\aut{\operatorname{Aut}}
\def\lcm{\operatorname{lcm}}
\def\sgn{\operatorname{sgn}}
\def\spn{\operatorname{span}}
\def\md{\operatorname{mod}}
\def\Norm{\operatorname{Norm}}
\def\dim{\operatorname{dim}}
\def\det{\operatorname{det}}
\def\Vol{\operatorname{Vol}}
\def\rk{\operatorname{rk}}
\def\ord{\operatorname{ord}}
\def\ker{\operatorname{ker}}
\def\div{\operatorname{div}}
\def\Gal{\operatorname{Gal}}
\def\GL{\operatorname{GL}}
\def\SNR{\operatorname{SNR}}
\def\WR{\operatorname{WR}}
\def\IWR{\operatorname{IWR}}
\def\scg{\operatorname{\left< \Gamma \right>}}
\def\swrh{\operatorname{Sim_{WR}(\Lambda_h)}}
\def\ch{\operatorname{C_h}}
\def\cht{\operatorname{C_h(\theta)}}
\def\scgt{\operatorname{\left< \Gamma_{\theta} \right>}}
\def\scgmn{\operatorname{\left< \Gamma_{m,n} \right>}}
\def\gat{\operatorname{\Omega_{\theta}}}
\def\mn{\operatorname{mn}}
\def\disc{\operatorname{disc}}
\def\rot{\operatorname{rot}}
\def\Prob{\operatorname{Prob}}
\def\co{\operatorname{co}}
\def\Ker{\operatorname{Ker}}

\newcommand{\cerr}[2]{{\color{green}#1} {\color{red}#2}}
\newcommand{\cimp}[2]{{\color{blue}#1} {\color{red}#2}}
\newcommand{\comm}[1]{\text{\ \color{cyan}#1}}

\section{Introduction}
\label{intro}

A classical arithmetic problem in the theory of finitely generated groups and rings is the study of the asymptotic order of growth of the number of subgroups of bounded index. A common approach to this problem involves studying the analytic properties of a corresponding zeta function (a Dirichlet series generating function) and then applying a Tauberian theorem to deduce information about the number in question, represented by the coefficients of this zeta function. This research direction received a great deal of attention over the years as can be seen from \cite{lubot}, \cite{sautoy_book}, \cite{sautoy1}, \cite{sautoy2}, \cite{reiner}, \cite{solomon} and the references within. In the recent years, a similar approach has also been applied to the more geometric setting of counting sublattices in lattices, e.g.~\cite{petrograd}, \cite{kuehnlein}, \cite{lf_pams}, \cite{scharlau}, \cite{shparlinski}.

In this note, we consider some special cases of the following general setting. Let $R$ be a commutative ring with identity such that for every natural number $n$ 
the set of ideals in $R$ of index $n$ is a finite number, call this number 
$a_n(R).$ One is often interested in the asymptotic behaviour of the sequence $a_n(R)$ or the summatory
sequence $A_N(R):= \sum_{n\leq N}a_n(R).$ To that end one can introduce the 
zeta function
$$\zeta(R,s) := \sum_n \frac{a_n(R)}{n^s}$$
and calculate its abscissa of convergence, so that one might apply a
Tauberian theorem.

We will make use of the following well-known result, which is a consequence of
Delange's extended version of Ikehara's Tauberian Theorem (\cite{delange}):

\begin{thm}\label{Tauber}
Let $(a_n)_{n\geq 1}$ be any sequence of non-negative real numbers and 
$$Z(s) := \sum_{n=1}^\infty \frac{a_n}{n^s}.$$
Assume that $Z(s)$ has abscissa of convergence $\sigma>0$ and admits a 
meromorphic extension to some neighborhood of the line $\Re(s)=\sigma$ with 
a pole of order $w$ at $s=\sigma$ and no other singularity. Then 
$$\sum_{n=1}^N a_n\sim c\cdot N^\sigma \cdot (\log N)^{w-1}$$
holds for $c= \frac{{\rm Res}({Z}, {\sigma})}{\sigma\cdot \Gamma(w)},$ where ${\rm Res}(Z, \sigma)$ stands for the residue of $Z(s)$ at $s=\sigma$ and $\Gamma(w)$ is the value of gamma function at $w$.
\end{thm}

The classical situation is that of rings of integers in number fields, where 
the zeta function $\zeta(R,s)$ is Dedekind's zeta function which converges 
for $\Re(s)>1$ and is meromorphic on $\mathbb C$ with a simple pole at $s=1.$ 
In this situation $A_N(R) \sim {c} \cdot N,$ the constant being 
the residue of $\zeta(R,s)$ at $s=1.$

One of our general results is the calculation of the zeta-function for 
polynomial rings:

\begin{thm} \label{Main Theorem}
Let $K$ be a number field with ring of
integers ${\O}_K.$ Then we have
$$\zeta({\O}_K[X] , s) = \prod_{d=1}^\infty \zeta_K(d(s-1)),$$
where $\zeta_K:= \zeta({\mathcal O}_K,\cdot)$ is the Dedekind zeta-function of 
the number field $K$.
\end{thm}

In particular, this function has abscissa of convergence $\sigma = 2$, 
meromorphic extension to $\Re(s)>1$ and simple poles at 
$s=\frac{d+1}d,\ d\in \Nn.$
As the poles accumulate towards 1, $\zeta(\O_K[X],s)$ cannot be extended 
meromorphically beyond $\Re(s) >1.$ The largest pole is at $s=2$
and therefore the number of ideals in $\O_K[X]$ of index less than $N$ grows like
a multiple of $N^2.$

We will, however, be more generally interested in ideals in quotient rings of 
${\O}_K[X],$ which are related to lattices in $\O_K^d$ for some $d.$ 
To that end, we will need some more machinery. Let $R$ {again} be any commutative ring.
If $I\leq R$ is an ideal of index $mn$ with coprime $m$ and $n$, then by the 
Chinese Remainder Theorem $R/I$ is a direct product of two rings of orders $m$ 
and $n$, respectively, and therefore $I$ is the intersection of two ideals of 
indices $m$ and $n$, respectively. These are uniquely determined by $I$, and 
therefore
$$a_{mn}(R) = a_m(R)\cdot a_n(R).$$
This means that the sequence $(a_n(R))_n$ is multiplicative, and so $\zeta(R,s)$ formally 
has a decomposition as an Euler product:
$$\zeta(R,s) = \prod_{p\in \mathbb P} E_p(R,s),\ \ \text{where}\ \ 
E_p(R,s) :=\sum_{k=0}^\infty \frac{a_{p^k}(R)}{p^{ks}}.$$
If $\tilde R\subseteq R$ is a subring of finite index, then for almost all prime
numbers the Euler-factors $E_p(R,s)$ and $E_p(\tilde R,s)$ coincide, and 
therefore the zeta functions $\zeta(R,s)$ and $\zeta(\tilde R,s)$ have the 
same convergence behaviour, if these finitely many Euler factors do not behave 
very abnormally. We will have to control this behaviour, when it comes to 
applications of this observation.

This paper is organized as follows. In Section~\ref{fixed}, we describe the construction of
lattices from ideals in rings of integers of number fields and quotient polynomial rings via the coefficient embedding.
We then use the Dedekind zeta function to obtain counting and density estimates on the numbers of ideal lattices
coming from a fixed ring of algebraic integers, improving on previous results~\cite{lind} in some cases. 
Our main result in this section (Theorem~\ref{sep_pol_count}) is an asymptotic estimate on the number of ideals in the quotient polynomial ring $\ZZ[X]/(f)$, where $f(X) \in \ZZ[X]$ is a monic separable polynomial. We also briefly comment on the situation when $f(X)$ is not separable. In Section~\ref{all_ideal}, we discuss the question of which sublattices of $\zed^d$ arise as ideal lattices from quotient polynomial rings $\zed[X]/(f)$ for some monic polynomial $f$ of degree~d. Specifically, we investigate the analytic properties of the corresponding zeta function (Theorem~\ref{zeta-d}) and, as a consequence of this theorem, show that the number of such ideal sublattices of index $\leq N$ grows asymptotically like $O(N^2)$ (Corollary~\ref{zeta-d-cor}). Finally in Section~\ref{zetaproof}, we prove Theorem~\ref{Main Theorem} and use it to deduce that the proportion of ideal lattices of index $\leq N$ among all sublattices of $\zed^d$ for $d \geq 3$ tends to 0 as $N \to \infty$ (Corollary~\ref{4-2}). 

\begin{rem} 
After finishing our calculations we became aware of the so-called 
K\"ahler zeta-functions, in particular the results {of} {Lustig in}{} 
\cite{lustig}. {With some additional work}, Lustig's results can be used to 
prove our Theorem \ref{Main Theorem}. 
His approach {however} is purely local, while we employ more global methods. 
Furthermore, 
our result is more specific {and neither} our Theorem \ref{zeta-d} {nor} 
Remark \ref{Cohen-Lenstra} (b) seem to have a parallel there. 
\end{rem}
\bigskip

\section{Ideal lattices from fixed rings}
\label{fixed}

In this section, we motivate our subsequent questions by looking at rings of 
integers in number fields.

Let $K$ be a number field of degree $d > 1$ with $r_1$ real and $r_2$ pairs of 
complex conjugate embeddings (so $d=r_1+2r_2$), and write $\O_K$ for its ring 
of integers. There exists $\theta \in \O_K$ such that $K=\que(\theta)$. Fix 
this $\theta$, then for each $\alpha \in K$ there exist $a_0,\dots,a_{d-1} \in 
\que$ such that
$$\alpha = \sum_{n=0}^{d-1} a_n \theta^n.$$
Notice that $\zed[\theta]$ is a finite-index subring of $\O_K$, possibly proper.
Define an embedding $\rho_K : K \to \que^d$ by 
$$\rho_K(\alpha) = (a_0,\dots,a_{d-1}),$$
then for every nonzero ideal $I \subseteq \O_K$ the image $\rho_K(I)$ is a 
full-rank 
lattice in $\real^d$, and for every ideal $I \subseteq \zed[\theta]$ the image 
$\rho_K(I)$ is a finite index sublattice of $\zed^d$. Furthermore, if some 
finite index sublattice $\Lambda \subseteq \zed^d$ is equal to $\rho_K(I)$ for 
some ideal $I \subseteq \O_K$, we will say that $\Lambda$ is an 
{\it $\O_K$-ideal lattice}, or just an ideal lattice when the choice of $K$ 
and $\theta$ is fixed. The first observation (see equation~(1) of~\cite{lind}) 
is that for each ideal $I \subseteq \O_K$, its norm is given by
$${\nu}(I) = {D_K^{-1}} \det(\rho_K(I)),$$
where $D_K :=  \det(\rho_K(\O_K))$.

For $T \in \real_{>0}$, define
\begin{eqnarray}
\label{NKT}
{\mathcal N}_K(T) & = & \left| \left\{ \Lambda \subseteq \zed^d : \Lambda = \rho_K(I) 
\text{ for some ideal } I \subseteq \O_K,\ \det \Lambda \leq T \right\} 
\right| \nonumber \\
& \leq & \left| \left\{ I \subseteq \O_K : \nu(I) \leq T{D_K^{-1}} \right\} \right|.
\end{eqnarray}
An upper bound on $\N_K(T)$ has been obtained in Theorem~2 of~\cite{lind}, 
which was then used in Corollary~1 of~\cite{lind} to establish a density 
estimate on ideal lattices from $\O_K$ among all sublattices of $\zed^d$. In 
the special case when $\O_K=\zed[\theta]$, one can easily use the standard 
analytic method to prove an asymptotic formula for $\N_K(T)$ as $T \to \infty$ 
and also deduce a more precise density estimate.

\begin{lem} \label{zeta} 
Suppose that $\O_K=\zed[\theta]$. As $T \to \infty$,
\begin{equation}
\label{NKT-bnd}
\N_K(T) \sim \frac{2^{r_1+r_2} \pi^{r_2} h_K R_K}{\omega_K \sqrt{|\Delta_K|}} T,
\end{equation} 
where $h_K$ is the class number, $R_K$ the regulator, $\Delta_K$ the 
discriminant, and $\omega_K$ the number of roots of unity of the number 
field~$K$.
\end{lem}

\proof
Note that $D_K=$ since $\O_K=\zed[\theta]$. Furthermore, every ideal lattice 
from $\O_K$ is a sublattice of $\zed^d$, and so there is equality in the 
inequality of~\eqref{NKT}. Let 
$s \in \cee$ and define the Dedekind zeta-function of $K$ by
$$\zeta_K(s) = \sum_{I \subseteq \O_K} {\nu}(I)^{-s} = \sum_{n=1}^{\infty} a_n n^{-s},$$
where $a_n$ is the number of ideals of norm $n$ in $\O_K$, and so, by 
\eqref{NKT},
$$\N_K(T) = \sum_{n=1}^{[T]} a_n.$$
It is well-known that $\zeta_K(s)$ is analytic in the half-plane $\Re(s) > 1$ 
and has only a simple pole at $s=1$ with residue given by the analytic class 
number formula:
$$\lim_{s \to 1} (s-1) \zeta_K(s) = \frac{2^{r_1+r_2} \pi^{r_2} h_K R_K}{\omega_K 
\sqrt{|\Delta_K|}}.$$
Combining this formula with the Tauberian Theorem~\ref{Tauber} yields the result.
\endproof

\begin{cor} \label{density} 
Define
$$\M(T) = \left| \left\{ \Lambda \subseteq \zed^d : \det \Lambda \leq T 
\right\} \right|.$$
Suppose that $\O_K=\zed[\theta]$. As $T \to \infty$,
$$\frac{\N_K(T)}{\M(T)} \sim \frac{2^{r_1+r_2} \pi^{r_2} h_K R_K  d}{\omega_K 
\sqrt{|\Delta_K|} \prod_{n=2}^d \zeta(n)} T^{1-d}.$$
\end{cor}

\proof
The zeta-function of all finite index sublattices of $\zed^d$ is
$$\zeta_{\zed^d}(s) = \sum_{\Lambda \subseteq \zed^d} (\det \Lambda)^{-s}.$$
It is a well-known fact (see, for instance p.~793 of~\cite{sautoy1}) that
$$\zeta_{\zed^d}(s) = \zeta(s) \zeta(s-1) \cdots \zeta(s-d+1),$$
where $\zeta(s)$ is the Riemann zeta-function. Hence $\zeta_{\zed^d}(s)$ is 
analytic in the half-plane $\Re(s) > d$ and has a simple pole at $s={d}$ with 
the residue
$$\zeta(d) \zeta(d-1) \cdots \zeta(2).$$
Combining this formula with Theorem~\ref{Tauber} implies that
$$\M(T) \sim \frac{\prod_{n=2}^d \zeta(n)}{d} T^d.$$
This estimate together with Lemma~\ref{zeta} yields the result.
\endproof

In the special case when $\O_K=\zed[\theta]$ there is also another way to look at 
ideal lattices from~$\O_K$. Let $f(X)$ be the minimal polynomial of $\theta$, 
which is a monic irreducible polynomial in $\zed[X]$ of degree~$d$. Then 
$\O_K \cong Z_f := \zed[X]/(f(X))$ and ideal lattices from $\O_K$ correspond to
ideal lattices in ~$Z_f$, which is a special case of the ideal 
lattice construction of~\cite{lub1}. Define a $\zed$-module isomorphism 
$\rho_f : Z_f \to \zed^d$, given by
$$\rho_f \left( \sum_{n=0}^{d-1} a_nx^n \right) = (a_0,\dots,a_{d-1}).$$
Let $I \subseteq Z_f$ be an ideal, then it is known that $\rho(I)$ is a finite 
index sublattice of $\rho(Z_f) = \zed^d$ (Lemma~3 of~\cite{lub1}). A counting 
estimate on such ideal sublattices is then given by our Lemma~\ref{zeta} above.

This leads us to the more general question: {if} $f(X) \in \zed[X]$ is a monic 
polynomial of degree $d$, what is the asymptotic behaviour of the number of
ideals of index at most $N$ in $Z_f:=\zed[X]/(f)$ as $N$ goes to infinity? We again define
$$\zeta(Z_f,s) = \sum_{I \subseteq Z_f} \nu(I)^{-s},$$
where the sum goes over finite index ideals and $\nu(I) := |Z_f/I|$. 

\begin{thm} \label{sep_pol_count} 
Let $f\in \zed[X]$ be monic and assume that 
{$f=g_1\cdots g_k$} with {$g_1,\dots,g_k$ in $\zed[X]$} irreducible, monic and 
pairwise distinct. Then the zeta-function $\zeta_f(s)$ converges for 
$\Re(s)>1,$ has a meromorphic extension to the halfplane
$${\{s\in \mathbb C \mid \Re(s) > 0\},}$$ 
and has a pole of order $k$ {at} $s=1.$ In particular,  
$A_N(Z_f) \sim c N(\log N)^{k-1}$ for some constant $c$.
\end{thm}

\proof
For $1\leq i \leq k,$ let $\theta_i\in \mathbb C$ be a root of $g_i $ and 
${\O}_i\subseteq \mathbb C$ be the integral closure of $\mathbb Z[\theta_i].$
We first observe that by a variant of the Chinese Remainder Theorem, the ring
$Z_f$ is isomorphic to a subring of finite index in 
$S:= {\O}_1 \times \cdots\times {\O}_k.$ If $d$ is the degree of $f,$ then both rings 
have an additive group isomorphic to $\zed^d$. Let $M$ be the index of $Z_f$ in $S.$
We denote by $a_n(Z_f)$ the number of ideals of index $n$ in $Z_f$ and by 
$a_n(S)$ the number of ideals of index $n$ in $S$ respectively. 

Whenever $I\subseteq S$ is an ideal of index $n$ in $S$, then $MI\subseteq Z_f$ 
is an ideal of index $M^{d-1}\cdot n$ in $Z_f.$ 
As multiplication by $M$ is injective on $S,$ this multiplication map from 
ideals in $S$ to ideals in $Z_f$ is injective, and for every natural number $N$ 
$$\sum_{n=1}^{NM^{d-1}} a_n(Z_f) \geq \sum_{n=1}^Na_n(S).$$
Therefore the abscissa of convergence of $\zeta(Z_f,s)$ is at least as 
large as that of $\zeta(S,s).$ This only depends on the fact, that the rings 
have an additive group isomorphic to $\mathbb Z^d$ and therefore holds for all 
subrings of finite index in $S.$ Furthermore, as remarked 
in Section~\ref{intro} the Euler-factors of both zeta-functions do coincide for 
all prime numbers not dividing the index.

We now consider Euler-factors of the ring $\zed + MS$, and show that they have the same 
convergence behaviour as Euler-factors for $S$. Since
$$\zed + MS \subseteq Z_f \subseteq S,$$
it would follow that Euler-factors of $Z_f$ must also have the same convergence
behaviour. Writing $M=p_1\cdot p_2\cdots p_l$ as a product of prime numbers, we see that 
$$S\supseteq \mathbb Z+p_1S \supseteq \mathbb Z+p_2(\mathbb Z+p_1S) = 
\zed + p_1p_2S\supseteq\cdots,$$
and this shows by the same argument as above that it suffices to treat the case $M=p$ and
the Euler-factor at the prime $p$.

Therefore let $S$ be a ring additively isomorphic to $\mathbb Z^d,$ $p$ a 
prime number and $R=\mathbb Z+pS.$ 
For every ideal $I$ of index $p^k$ in $R$, $SI$ is an ideal of $p$-power index 
in $S.$ As $pS \subseteq R$, we see that $pSI \subseteq I$ and that the index of 
$SI$ in $S$ is at least $p^{-d}\cdot p^k.$ For two ideals $I_1,I_2$ in $R$ the 
equality
$SI_1 = SI_2$ implies $pI_2\subseteq pSI_2 = pSI_1 \subseteq I_1$ and by symmetry
$pI_1\subseteq I_2.$ Hence $p^2I_1 \subseteq pI_2 \subseteq I_1.$ The 
cardinality of the fibre of the map $I\mapsto SI$ therefore is at most the 
number of subgroups in $I_1/p^2I_1 \cong \mathbb Z^d / p^2\mathbb Z^d$, which is 
finite and independent of $I_1.$ Call this number $\gamma$. 

Then 
$${\sum_{k=0}^K a_{p^k}(S) \geq \frac1\gamma \sum_{k=0}^{K+d} a_{p^k}(\mathbb Z+pS),}$$
thereby showing that the Euler-factor in $\zeta(\zed+pS,s)$ converges whenever
that of $\zeta(S,s)$ does. 
Coming back to our statement, the Euler-factors in the zeta-function for 
$\O_1\times \dots\times \O_k$ are just the products of the Euler-factors 
of $\zeta(\O_1,s),\dots , \zeta(\O_k,s),$ which all converge for $\Re(s)>0.$
This completes the proof.
\endproof

\begin{rem} 
What if $f(X)$ is not separable? Here are some speculations. 
Going through the proof of Theorem 2.3, it is sufficient to understand the 
zeta-functions for rings of the type $\zed [X] / (f^e),$ where $f$ is monic 
and irreducible. We calculated this for ${f(X)=X}, e=2$ and obtained
$$\zeta(\zed[X]/(X^2) ,s) = \zeta(s)\cdot \zeta(2s-1),$$
which again has a double pole at $s=1.$ We expect 
$$\zeta(\zed[X]/(X^e) ,s) = 
\zeta(s)\cdot \zeta(2s-1)\cdot \zeta(3s-2)\cdots \zeta(es-(e-1)),$$
which has a pole of order $e$ at $s=1.$ More generally, we would expect that the order of 
the pole of 
$\zeta(\zed[X]/(f),s)$ at $s=1$ is the sum of the multiplicities of the 
irreducible monic factors of $f.$
\end{rem}
\bigskip

\section{Ideal lattices in dimension $d$}
\label{all_ideal}

In this section, we aim to understand which subgroups of $\zed^d$ arise as 
images of finite index ideals in $Z_f:=\zed[X]/(f)$ under $\rho_f$ for some 
monic polynomial $f(X)$ of degree $d.$ Again, we fix
the group isomorphism $\rho_f : Z_f \to \zed^d$, given by
$$\rho_f \left( \sum_{n=0}^{d-1} a_nX^n \right) = (a_0,\dots,a_{d-1}).$$
Like every finite index subgroup of $\zed^d$, such a subgroup is generated by 
the columns of an integral full-rank upper triangular matrix 
$$A:= \begin{pmatrix} a_{0,0} & a_{0,1} & \dots & a_{0,d-1} \\
                        0    & a_{1,1} & \dots & a_{1,d-1} \\
                        \vdots & \ddots & \ddots & \vdots \\
                        0 & \dots & 0 & a_{d-1,d-1} \end{pmatrix},$$
the $j$-th column corresponding to a generator of degree $j$ in the ideal.

Multiplication by $X$ is an endomorphism of the ideal. This shows that for each 
$0\leq j\leq d-2$ the $j$-th column of $A$ shifted down by one belongs to the 
subgroup generated by the $0$-th to $(j+1)$-th column. For $j=d-1$ it shows, 
that there exists some monic poynomial $f$
of degree $d$ such that the $(d-1)$-th column shifted by one (in 
$\zed^{d+1}$) 
belongs to the subgroup generated by the columns of $A$ and the column 
containing the coefficients of $f.$ This last demand is always satisfied, and 
we therefore only have to deal with the first set of conditions!

We call such a matrix an {\it  idealizing matrix}, and two idealizing matrices 
are called {\it indifferent}, whenever they correspond to the same ideal. 
One can check by induction on $d$ that an idealizing upper triangular matrix satisfies: 
\begin{equation} \label{Divisibility}
\forall 1\leq j\leq d-1: a_{j,j}\ \text{divides}\ a_{i,k}\ \text{whenever}\ 
0\leq i,k \leq j.
\end{equation}
Details of this calculation can be found in \cite{schwerdt}.
The idealizing matrices for which in every $i$-th row the entries
are between $0$ and $a_{i,i}-1$ give a set of representatives for the 
indifferency-classes. We call these matrices {\it reduced} idealizing matrices,
and now count them by induction in the following way.

\begin{lem} \label{New Columns}
Let $A\in \mathbb \zed^{d\times d}$ be a reduced idealizing matrix. Then the 
number of $\beta\in \zed ^d$ such that 
$$\tilde A {\ :=} \begin{pmatrix} A & \beta \\0 & 1 \end{pmatrix} \in\zed^{d+1\times d+1}$$
is a reduced idealizing matrix is $a_{d-1,d-1}^{d}.$
\end{lem}
\proof

As the first $d-1$ columns of $\tilde A$ already satisfy the conditions imposed 
on idealizing matrices, we only have to care about the condition that the 
shifted $d$-th column is contained in the $\zed$-span of the columns of
$\left(\begin{smallmatrix}\delta \\ a \end{smallmatrix}\right)$ 
with $\delta\in \zed^{d-1}, a=a_{d-1,d-1}\in \zed.$ 

Then our condition is
$$\begin{pmatrix} 0\\ \delta\end{pmatrix} - a \beta \in A\zed^d.$$
Note that by \eqref{Divisibility} $a$ divides every entry of $A$ and that 
therefore we get
$$\beta \in \frac1a \left(\begin{pmatrix} 0\\ \delta\end{pmatrix} + A\zed^d\right).$$
Two different values of $\beta$ define the same indifferency class if and only
if their difference is in $A\zed^d.$ As the index of $A\zed ^d$ in 
$\frac1aA \zed^d$ is $a^d,$ we get exactly $a^d$ possible indifferency classes of
matrices $\tilde A$. \endproof

\begin{thm} \label{zeta-d} 
Let $d\in \mathbb N$ be some natural number, and denote by
$c_n^{(d)}$ the number of subgroups in $\zed^d$ of index $n$ which are in the image of 
$\rho_f$ for some monic polynomial $f\in \zed[X]$ of degree $d.$
Then 
$$\zeta^{(d)}(s) := \sum_{n=0}^\infty \frac{c_n^{(d)}}{n^s} = 
\left[\prod_{i=1}^{d-1} \zeta(i(s-1))\right] \cdot \zeta(ds),$$
{where} $\zeta = \zeta^{(1)}$ is the Riemann zeta-function.
\end{thm}

\proof 
The proof is by induction on $d.$ We have to count indifferency classes 
of matrices with given determinant $n.$ The diagonal entries of such a matrix 
are numbers $a_0,\dots ,a_{d-1},$ where each $a_i$ divides the $a_j$ with $j<i.$
It is more convenient to write these numbers as products
$${b_1 b_2\cdots b_d , \ b_2 \cdots b_d,\ \dots,\ b_{d-1}b_d,\ b_d,}$$
where $b_1,\dots, b_d$ are any natural numbers.

For $d=1,$ we just have one indifferency class of matrices of determinant $n$ 
for every $n,$ hence
$$\zeta^{(1)}(s) = \zeta(s).$$
For $d=2$ and every value of $b_1$, by Lemma \ref{New Columns}, we obtain 
exactly $b_1$ classes of matrices
$$\begin{pmatrix} b_1 & \beta \\ 0 & 1 \end{pmatrix},$$
and every matrix with arbitrary $b_1,b_2$ is $b_2$ times one of these, therefore
we have $b_1$ matrices with fixed $b_1.$ 
This gives
$$\zeta^{(2)}(s) = \sum_{b_1,b_2} \frac{b_1}{b_1^s\cdot b_2^{2s}} = 
\sum_{b_1} \frac1{b_1^{s-1}} \cdot \sum_{b_2} \frac1{b_2^{2s}} = 
\zeta(s-1) \times \zeta(2s).$$
Going on like this, using Lemma \ref{New Columns} repeatedly, we find that for 
fixed $b_1,\dots ,b_d$ we have $b_1\cdot b_2^2 \cdot b_3^3 \cdot b_{d-1}^{d-1}$ 
indifferency classes of idealizing matrices with diagonal
$${(b_1\cdots b_d,\ b_2\cdots b_d,\ \dots,\ b_{d-1}b_d,\ b_d)},$$
and this leads to 
$$\zeta^{(d)}(s) = \sum_{b_1,\ldots ,b_d = 1}^\infty
\frac {b_1\cdot b_2^2 \cdot b_3^3 \cdots  b_{d-1}^{d-1}}{(b_1\cdot b_2^2 \cdot b_3^3 \cdots b_{d-1}^{d-1}\cdot b_d^d)^s} =
\prod_{i=1}^{d-1} \zeta(i(s-1)) \times 
\zeta(ds).$$
\endproof

\begin{cor} \label{zeta-d-cor} 
For $d\geq 2,$ the number $A_N$ of subgroups in $\mathbb Z^d$ of
index less than or equal to $N$, which are in the image of 
$\rho_f$ for some monic polynomial $f\in \zed[X]$ of degree $d$, grows 
asymptotically as $cN^2,$ where 
$$c = \frac{{\rm Res}(\zeta^{(d)},2)}2 = \left[
\prod_{i=2}^{d-1} \zeta(i)\right] \cdot \zeta(2d).$$
\end{cor}

\proof
Apply the Tauberian Theorem \ref{Tauber} to $\zeta^{(d)}.$ The abscissa of convergence 
is $2$ and $\Gamma(2) = 1$; use the fact that the Riemann zeta-function $\zeta(s)$ has 
a simple pole at $s=1$ and no other poles in $\mathbb C.$
\endproof
\bigskip

\section{Zeta-functions of Polynomial Rings}
\label{zetaproof}

In this section we prove Theorem~\ref{Main Theorem}. We start with its version for $K=\que$, and then explain how to generalize it to any number field.

\begin{prop} 
We have 
$$\zeta(\zed[X],s) = \prod_{i=1}^\infty \zeta(i(s-1)).$$
\end{prop}

\proof 
Every ideal of finite index in $\zed[X]$ contains some monic polynomial,
and therefore $a_n(R) = \lim_{d\to\infty} a_n^d.$ The discussion in the previous
section shows that for a fixed $n$ the sequence $(a_n^d)_d$ becomes constant for 
$d>n$ (very rough estimate). This leads to 
$$\zeta(\zed[X],s) = \lim_{d\to \infty} \zeta^{(d)}(s),$$
which gives the desired result. The infinite product on the right hand side 
in fact converges absolutely for $\Re(s) > 2.$
\endproof

\begin{cor} \label{4-2} 
For $d\geq 3,$ the proportion of the ideal lattices of index 
$\leq N$ in $\zed^d$ among all subgroups of index $\leq N$ tends to zero as 
$N$ goes to infinity.
\end{cor}
\proof The zeta-function $\zeta(\zed[X] ,s)$ converges for $\Re(s) >2$ and has
a simple pole at $s=2.$  Theorem \ref{Tauber} says that the number of ideal 
lattices has quadratic growth, while that for all subgroups has growth of 
order $N^d.$
\endproof

\proof[Proof of Theorem~\ref{Main Theorem}]
We now explain how the calculations in the proof of the last Proposition can be extended to
rings of integers in arbitrary number fields instead of $\zed.$

\begin{itemize}
\item [(a)] 
Let $R$ be an infinite principal ideal domain such that every nonzero 
ideal in $R$ has finite index. Here, the method using idealizing matrices can 
be extended directly. Repeating {\it mutatis mutandis} the arguments from the 
last section shows the formal identity
$$\zeta(R[X],s) = \prod_{i=1}^\infty \zeta(R, i(s-1)).$$
\item [(b)]
If $\O_K$ is the ring of integers in a fixed number field $K$ of finite degree 
over $\mathbb Q,$ we use the Euler-product
$$\zeta(\O_K[X],s) = \prod_{p\in \mathbb P} \zeta(\zed_{(p)}\otimes_\zed\O_K[X] , s),$$
as the ideals of $p$-power index in $\O_K[X]$ correspond bijectively and index 
preservingly to ideals of finite index in $\zed_{(p)}\otimes_\zed{\O_K}.$ But this 
last
ring is a Dedekind domain with finitely many prime ideals only (i.e. semilocal) 
and therefore is a principal ideal domain which satisfies the condition from 
part (a) of this proof. Hence
$$\zeta(\O_K[X],s) = {\prod_{p\in\mathbb P}}\prod_{i=1}^\infty \zeta(\zed_{(p)}\otimes_\zed\O_K, i(s-1)).$$
For every prime number $p,$ the factor $\zeta(\zed_{(p)}\otimes_\zed\O_K, i(s-1))$
is the pro\-duct of the ``true'' Euler factors of $\zeta(\O_K, i(s-1))$ for 
prime ideals in $\O_K$ containing $p,$ and therefore, using the 
Euler-decomposition of
the Dedekind zeta-function $\zeta(\O_K,\cdot),$ we arrive at the assertion of
Theorem \ref{Main Theorem}.
\end{itemize}
\endproof

\begin{rem} \label{Cohen-Lenstra} 
We conclude this paper with some remarks.
\begin{itemize}
\item[(a)] 
The first factor in $\zeta(\O_K[X],s)$ is $\zeta(\O_K,s-1),$ which is the
Hasse-Weil zeta-function of the affine line {$\aaa^1_{\O_K}$}. 

This corresponds to those ideals in $\O_K[X]$ where the quotient ring has 
squarefree characteristic. It would be interesting to understand the 
zeta-function we have calculated from this point of view and perhaps describe it
as that of some deformation of the affine line.

There is another way to discover $\zeta(\O_K,s-1)$ as a factor in 
$\zeta(\O_K[X],s),$ namely by counting only ideals which contain a linear monic 
polynomial. 

\item[(b)]
If $a_n({\mathfrak G})$ denotes the number of isomorphism classes of abelian
groups of order $n$, then an argument based on the structure theorem for 
finitely generated abelian groups shows, that
$$\sum_{n=1}^\infty \frac{a_n({\mathfrak G})}{n^{s}} = 
\prod_{d=1}^\infty \zeta(ds) = 
\zeta(\zed[X],s+1).$$
Namely, for every choice $b_1,\dots ,b_d\in \mathbb N$ there is exactly one 
abelian group with elementary divisors 
$b_d, b_db_{d-1} , b_db_{d-1}b_{d-2}, \dots, b_d{\cdot \ldots \cdot} b_1.$
This group has order $b_1\cdot b_2^2 {\cdot \ldots \cdot} b_d^d,$ and therefore
we get 

\begin{align*}\sum_{n=1}^\infty \frac{a_n({\mathfrak G})}{n^{s}} & = \lim_{d\to\infty} 
\sum_{b_1,\dots ,b_d} \frac1{(b_1\cdot b_2^2\cdots b_d^d)^s}\\
& = \lim_{d\to\infty} \zeta(s)\cdot\zeta(2s)\cdots \zeta(ds).\end{align*}

This is closely related to the Cohen-Lenstra heuristics, cf.\ \cite{cohen}, 
although there the product $\prod_{i\geq 1}\zeta(s+i)$ plays a more important 
role.
However, the residue of our $\zeta(\O_K[X],s)$ at $s=2$ is the number 
$C_\infty$ on page 35 of loc. cit.
We thank Ernst-Ulrich Gekeler for kindly guiding us towards \cite{cohen}.

Maybe it is possible by such means to obtain a nice interpretation of $\zeta(\zed[X],s)$ 
as the zeta-function counting finite abelian groups endowed with some 
endomorphism satisfying certain properties. 
The number $a_n(\mathfrak G)$ is the number of orbits of subgroups of index $n$ 
in $\zed^\infty = \bigoplus_{i\in\mathbb N} \zed$ under the action of the 
automorphism group of $\zed^\infty.$ 
It would also be interesting to count the number of orbits of 
$\text{Aut}(\zed[X])$ acting on ideals of index $n.$ We did not yet carry out 
the corresponding calculations.

\item[(c)] We should mention the short note \cite{witt} of Witt. His results 
 imply that the sum formally defining the $\zeta$-function for the 
polynomial ring $\mathbb Z[X,Y]$ does not converge anywhere in the complex 
plane.
\end{itemize}
\end{rem}


\end{document}